\documentclass[12pt]{article}
\author{F. Argentieri\footnote{FERNANDO ARGENTIERI fernando.argentieri@math.uzh.ch: Institut f\"ur Mathematik Universit\"at Z\"urich, Winterthurerstrasse 190, 8057 Zurich, Switzerland.}, B. Fayad\footnote{BASSAM FAYAD bassam@umd.edu: Department of Mathematics, University of Maryland, 4176
Campus Drive, College Park, MD 20742–4015, USA}}
\usepackage{amssymb,amsmath,amsfonts,bbm,pifont,upgreek,bbold,accents}

\newcommand{\norm}[1]{\bigl\| #1 \bigr\|}
\renewcommand{\norm}[1]{|\!|\!| #1 |\!|\!|}
 \def\paralll{|\!|}
 \renewcommand{\parallel}{\paralll}
\usepackage[normalem]{ulem}

\UseRawInputEncoding
\font\teneufm=eufm10
\font\seveneufm=eufm7
\font\fiveeufm=eufm5
\newfam\eufmfam
\textfont\eufmfam=\teneufm
\scriptfont\eufmfam=\seveneufm
\scriptscriptfont\eufmfam=\fiveeufm

\newcommand\beq[1]{ \begin{equation}\label{#1} }
\newcommand{\eeq}{ \end{equation} }

\newcommand{\beqa}{ \begin{align*} }
\newcommand{\eeqa}{ \end{align*} }

\newcommand{\beqano}{ \begin{eqnarray*} }
\newcommand{\eeqano}{ \end{eqnarray*} }

\newtheorem{theorem}{Theorem}
\newtheorem{definition}{Definition}
\newtheorem{proposition}{Proposition}
\newtheorem{lemma}{Lemma}
\newtheorem{sublemma}{Sublemma}
\newtheorem{remark}{Remark}
\newtheorem{notationalremark}{Notations}
\newtheorem{corollary}{Corollary}
\newtheorem{assumption}{Assumption}
\newtheorem{claim}{Claim}
\newtheorem{tools}{$\negsp\negsp$}[subsection]

\newcommand\thm[1]{ \begin{theorem}\label{#1}}
\newcommand\thmtwo[2]{ \begin{theorem}[#1]\label{#2}}
\newcommand\ethm{ \end{theorem} }
\newcommand\dfn[1]{ \begin{definition}\label{#1} \rm}
\newcommand\dfntwo[2]{ \begin{definition}[#1]\label{#2} \rm}
\newcommand\edfn{ \end{definition} }
\newcommand\pro[1]{ \begin{proposition}\label{#1}}
\newcommand\protwo[2]{ \begin{proposition}[#1]\label{#2}}
\newcommand\epro{ \end{proposition} }
\newcommand\lem[1]{ \begin{lemma}\label{#1}}
\newcommand\lemtwo[2]{ \begin{lemma}[#1]\label{#2}}
\newcommand\elem{ \end{lemma} }
\newcommand\sublem[1]{ \begin{sublemma}\label{#1}}
\newcommand\sublemtwo[2]{ \begin{sublemma}[#1]\label{#2}}
\newcommand\esublem{ \end{sublemma} }
\newcommand\rem[1]{ \begin{remark}\label{#1} \rm}
\newcommand\erem{ \end{remark} }
\newcommand\notrem[1]{ \begin{notationalremark}\label{#1} \rm}
\newcommand\enotrem{ \end{notationalremark} }
\newcommand\cor[1]{ \begin{corollary}\label{#1}}
\newcommand\cortwo[2]{ \begin{corollary}[#1]\label{#2}}
\newcommand\ecor{ \end{corollary} }
\newcommand\asmp[1]{ \begin{assumption}\label{#1}}
\newcommand\asmptwo[2]{ \begin{assumption}[#1]\label{#2}}
\newcommand\easmp{ \end{assumption} }
\newcommand\clm[1]{ \begin{claim}\label{#1}}
\newcommand\eclm{ \end{claim} }
\newcommand{\proof}{\par\medskip\noindent{\bf Proof. }}
\expandafter\chardef\csname pre amssym.def
at\endcsname=\the\catcode`\@
\catcode`\@=11
\def\undefine#1{\let#1\undefined}
\def\newsymbol#1#2#3#4#5{\let\next@\relax
 \ifnum#2=\@ne\let\next@\msafam@\else
 \ifnum#2=\tw@\let\next@\msbfam@\fi\fi
 \mathchardef#1="#3\next@#4#5}
\def\mathhexbox@#1#2#3{\relax
 \ifmmode\mathpalette{}{\m@th\mathchar"#1#2#3}%
 \else\leavevmode\hbox{$\m@th\mathchar"#1#2#3$}\fi}
\def\hexnumber@#1{\ifcase#1 0\or 1\or 2\or 3\or 4\or 5\or 6\or 7\or
8\or
 9\or A\or B\or C\or D\or E\or F\fi}
\ifcase\@ptsize
 \font\tenmsb=msbm10
 \font\sevenmsb=msbm7
 \font\fivemsb=msbm5
\or
 \font\tenmsb=msbm10 scaled \magstephalf
 \font\sevenmsb=msbm7 scaled \magstephalf
 \font\fivemsb=msbm5  scaled \magstephalf
\or
 \font\tenmsb=msbm10 scaled \magstep1
 \font\sevenmsb=msbm7 scaled \magstep1
 \font\fivemsb=msbm5 scaled \magstep1
\fi
\newfam\msbfam
\textfont\msbfam=\tenmsb
\scriptfont\msbfam=\sevenmsb
\scriptscriptfont\msbfam=\fivemsb
\edef\msbfam@{\hexnumber@\msbfam}
\def\Bbb#1{\fam\msbfam\relax#1}
\def\widehat#1{\setboxz@h{$\m@th#1$}%
 \ifdim\wdz@>\tw@ em\mathaccent"0\msbfam@5B{#1}%
 \else\mathaccent"0362{#1}\fi}
\def\widetilde#1{\setboxz@h{$\m@th#1$}%
 \ifdim\wdz@>\tw@ em\mathaccent"0\msbfam@5D{#1}%
 \else\mathaccent"0365{#1}\fi}

\def\RIfM@{\relax\ifmmode}
\def\nonmatherr@#1{\errmessage{\string#1\space allowed only in math mode}}
\def\Bbb{\RIfM@\expandafter\Bbb@\else
 \expandafter\nonmatherr@\expandafter\Bbb\fi}
\def\Bbb@#1{{\Bbb@@{#1}}}
\def\Bbb@@#1{\fam\msbfam\relax#1}
\def\setboxz@h{\setbox\z@\hbox}
\def\wdz@{\wd\z@}
\catcode`\@=\csname pre amssym.def at\endcsname
\newcommand{\giu}{{\medskip\noindent}}
\newcommand{\Giu}{{\bigskip\noindent}}
\newcommand{\nl}{{\smallskip\noindent}}

\newcommand{\qed}{\hskip.5truecm
\vrule width 1.7truemm height 3.5truemm depth 0.truemm
\par\Giu}

\newcommand{\negsp}{\hspace{-.09truecm}}  

\newcommand\su[1]{ \frac{1}{ {#1}} }

\renewcommand{\a }{ {\alpha}   }
\renewcommand{\b}{ {\beta}   }

\newcommand{\e }{ {\epsilon}   }

\def\R{\mathbb R}

\def\Z{\mathbb Z}
\renewcommand\subset{\subseteq}

\usepackage{color}

\title{Reducibility without KAM}
\begin{document}
\maketitle
\begin{abstract}
    We prove rotations-reducibility for close to constant quasi-periodic $SL(2,\Bbb{R})$ cocycles in one frequency in the finite regularity and smooth cases, and derive some applications to quasi-periodic Schr\"odinger operators. 
\end{abstract}
\subsection*{Introduction}
\nl
In this paper we will study smooth quasi-periodic $SL(2,\Bbb{R})$ cocycles in one frequency. These are skew products of the form
\beqano
(\a,A)&:&\Bbb{T}\times \Bbb{R}^2\rightarrow \Bbb{T}\times \Bbb{R}^2\\
&&(x,y)\rightarrow(x+\a, A(x)y),
\eeqano
$A \in C^\infty(\Bbb{T},SL(2,\Bbb{R}))$. 

We will be interested in the case $A$ is close to a constant matrix. A classical problem is to see if $(\a,A)$ is reducible, that is, conjugated to a constant matrix cocycle. We say that $(\a,A)$ is $C^\infty$-{\it reducible} if there exists $B \in C^\infty(\Bbb{R}/2\Bbb{Z}, SL(2,\Bbb{R}))$ and $A_* \in SL(2,\Bbb{R})$ such that 
\begin{eqnarray*}
{B(\theta+\a)} A(\theta) {B(\theta)}^{-1} = A_*, \qquad \forall \theta \in \Bbb{R}/2\Bbb{Z}.
\end{eqnarray*}
 
Reducibility is an important question in the study of quasi-periodic cocycles and its applications to the spectral theory of Schr\"odinger operators. 

When $\a$ satisfies a Diophantine condition, many reducibility results were obtained by the KAM (Kolmogorov-Arnold-Moser) technique since the seminal paper by Dinaburg and Sinai \cite{1}. Most of the results were obtained for real analytic cocycles, but KAM methods also yield  reducibility results in the smooth category, see \cite[Section 2.4]{2} for example. 

If $\a$ is just irrational, reducibility does not hold in general. Indeed, the rotation-valued cocycle $(\a,A)$,  $A(\cdot)=R_{\phi(\cdot)}$, with $\phi \in C^\infty(\Bbb{T},\Bbb{R})$ and
\beqano
R_{\phi(x)}:=
\begin{pmatrix}
\cos (2\pi\phi(x))  & -\sin (2\pi\phi(x)) \\
\sin (2\pi\phi(x)) & \cos (2\pi\phi(x))
\end{pmatrix},
\eeqano
is smoothly reducible if and only if the cohomological equation 
\begin{equation}Ê\label{eqconj} \tag{$\mathcal{E}$} \phi(x)-\int_{\Bbb{T}} \phi(\theta)d\theta=h(x+\a)-h(x)\end{equation} 
has a smooth solution. It is known that the cohomological equation does not have smooth solutions in general when $\a$ is not Diophantine.  

However, one can still ask about reducibility to a rotation-valued cocycle, that is, the existence of a smooth conjugacy $B$ such that
$B(x+\a)A(x)B^{-1}(x) \in SO(2,\Bbb{R})$. We then say that the cocycle is rotations-reducible. Rotations-reducibility is important in the global theory of 1 dimensional quasi-periodic cocycles, and  has many  applications for quasi-periodic Schr\"odinger operators, some of which will be mentioned below.

Rotations-reducibility for close to constant cocycles, irrespective of any arithmetic condition on the irrational base frequency $\a \in \Bbb{T}$, was obtained in \cite{3} in the real analytic category.  The aim of this paper is to extend the results of \cite{3} to the finite regularity and smooth case.  Rotations-reducibility was obtained in some ultradifferentiable classes (for example Gevrey) by  Cheng, Ge, You,  Zhou in \cite{4}, while our rotations-reducibility result holds in finite regularity (see also \cite{12} for questions about reducibility with low regularity).

Before stating our main theorem, we recall the definition of  the fibered rotation number $\rho=\rho(\a,A)$. If $A$ is homotopic to the identity, $G$ is the projective cocycle associated with $(\a,A)$, that is:
\beqano
G&:&\Bbb{T}\times \Bbb{S}^1\rightarrow\Bbb{T}\times \Bbb{S}^1\\
&&(x,y)\rightarrow(x+\a, \frac{A(x)y}{\parallel A(x)y\parallel}),
\eeqano
and $\tilde{G}(\theta, y)=(\theta+\a,y+f(\theta,y))$ is a lift of $G$ in $\Bbb{T}\times\Bbb{R}$, then:
\beqano
\lim_{n\rightarrow\infty}\su{n}\sum_{i=0}^{n-1}f(\tilde{G}^{i}(\theta,y))
\eeqano
exists and it is independent of $(\theta,y)$. The class of this number in $\Bbb{T}$ is the fibered rotation number.

\giu
\nl
We will prove the following
\thm{t1} \label{theorem1}
Let $\e>0, r_0\geq 200$. There exists $ \e_0>0$ such that, for $\a\in\Bbb{R}-\Bbb{Q}$ there exists a set $Q(\a)\subset\Bbb{T}$ with $m (Q(\a)^{c})<\e$ such that the following holds:

\giu
if $A_0\in C^{\infty}(\Bbb{T}, SL(2,\Bbb{R}))$, $R\in SO(2,\Bbb{R})$, 
$\rho(\a,A_0)\in Q(\a)$ and:
\beqano
\parallel A_0-R\parallel_{50r_0}<\e_0,
\eeqano
then, there exist $B\in C^{\infty}(\Bbb{T},SL(2,\Bbb{R}))$, $\phi\in C^{\infty}(\Bbb{T},\Bbb{R})$ such that,
\beqano
\label{eq.conj}ï¿½B(x+\a)A_0(x)B^{-1}(x)=R_{\phi(x)},
\eeqano
with:
\beqano \label{eq.small} 
\parallel\phi-\rho\parallel_{r_0}, \parallel B-Id.\parallel_{r_0}<\sqrt{\e_0}.
\eeqano
\ethm{}

\medskip

\noindent{\bf Finite regularity.} Observe that in the proof of Theorem \ref{theorem1}, in order to show the convergence of the scheme in norm $C^{r_0}$ we will use only derivatives of $A$ up to order $50r_0$. In particular, the same statement holds in class $C^{r_0}$ if $A$ is only assumed to be in $C^{50r_0}(\Bbb{T}, SL(2,\Bbb{R}))$. Our assumption that $r_0\geq 200$ is not optimal, neither is the closeness condition in class $C^{50r_0}$. In fact, we will explain later why the more $\a$ is Liouville, the  more $r_0$ can be taken small. 

\medskip 

\noindent{\bf Reducibility in the Diophantine case.}ÊWhen the base frequency $\a$ is Diophantine, our result implies smooth reducibility since the cohomological equation \eqref{eqconj} has a smooth solution $h$. Finite regularity reducibility results, depending on the Diophantine exponent of $\a$, can also be obtained from the proof of Theorem \ref{theorem1} and the study of \eqref{eqconj} in finite regularity.

\medskip 

\noindent {\bf Full measure rotations-reducibility. Eliasson's theory.}ï¿½ Eliasson developed a non-standard KAM scheme
which enabled him to prove a much stronger version of Dinaburg-Sinai theorem. He
showed that for every Diophantine $\a$, there exists a full Lebesgue measure set $Q(\a)$ such that if the cocycle is sufficiently
close to constant (depending on $\a$) and if $\rho(\a,A)\in Q(\a)$ then $(\a,A)$ is reducible \cite{eliasson}.

Note that for arbitrary irrational frequency $\a$, unconditional almost reducibility and therefore 
 full measure rotations-reducibility, for analytic close to constant cocycles, follows from Avila's global theory \cite{avila_global}.

In \cite{houyou}, Hou and You proved a continuous time version of \cite{3} with a different method. They also prove that almost reducibility always holds in the close to constants regime, which gives rotations reducibility for a full measure set of fibered rotation numbers. 

In the smooth category, full measure rotations-reducibility does not hold, since counter-examples were found in Gevrey class by Avila and Krikorian \cite{AKsmooth}. In this respect, the result of Theorem \ref{theorem1} is optimal.

\medskip

\noindent{\bf The case of general matrix cocycles.} 
Reducibility results for quasi-periodic cocycles valued in $GL(d,\mathbb{C})$ above a one-dimensional Liouvillean rotations, extending the results of \cite{3,houyou}, were obtained in \cite{WangZhou} in the real analytic setting. The scheme developed in the proof of Theorem \ref{theorem1} should be useful to address similar extensions in the smooth setting.

\bigskip 

\noindent{\bf Applications to one-dimensional quasi-periodic Schr\"odinger operators.} The main source of examples of cocycles that we are considering are the Schr\"odinger cocycles 
\beqano
A_{v,E}(x)=
\begin{pmatrix}
E-v(x)  & -1 \\
1 & 0
\end{pmatrix},
\eeqano
where $v \in C^\infty(\Bbb{T},\Bbb{R})$ and $E \in \R$ are related to the spectral study of one-dimensional quasi-periodic Schr\"odinger operators:
\beqano
(Hu)_n=u_{n+1}+u_{n-1}+v(\theta+n\a)u_n.
\eeqano

An application of our main result is the following. 
 \thm{t3} \label{theorem3}
 If $v\in C^{\infty}(\Bbb{T},\Bbb{R})$ is close to a constant, $\a\in\Bbb{R}$, then there is a positive measure set of $E$ such that $(\a,A_{v,E})$ is rotations-reducible. In particular, for $v$ close to a constant there exists some absolutely continuous part in the spectrum of the corresponding Schr\"odinger operator.
 \ethm

Another application concerns Schr\"odinger's conjecture stating that for general discrete Schr\"odinger operators over uniquely ergodic base dynamics, all generalized eigenfunctions are bounded for almost every energy in the support of the absolutely continuous part of the spectral measure. This conjecture was disproved by Avila (see \cite{5}). In \cite{17}, Marx and Jitomirskaya pointed out that the Schr\"odinger's conjecture in the quasiperiodic setting is still an open problem in the smooth category. The following result, implies that the Schr\"odinger conjecture is true for the one-dimensional smooth quasi-periodic Schr\"odinger operators over irrational circle rotations. 
   
\thm{t2} \label{theorem2} 
Let $v\in C^{\infty}(\Bbb{T}, \Bbb{R})$. Then, for almost every $E\in \Bbb{R}$, $(\a, A_{v,E})$ is either non uniformly hyperbolic or smoothly  rotations-reducible.
 \ethm
 The implication of Theorems \ref{theorem3} and  \ref{theorem2} from  Theorem \ref{theorem1} can be obtained in exactly the same way as  Theorem 1.1 is obtained from Theorem 1.3 in \cite{3}. 
 
\bigskip

\noindent{\bf Main novelties and outline of the proof.} The proof follows an inductive conjugation scheme based on the so-called cheap trick introduced in \cite{2} and further developed in \cite{3}.

One main novelty in our approach is that in finite regularity it is possible to fully use the strength of the cheap trick in which no loss of derivatives is incurred in finding the conjugacy (see \eqref{intro_eq1} below), since no cohomological equation is solved. As a consequence, the convergence of our conjugation scheme is quite different from that of KAM smooth schemes, and is much simpler. In particular, it does not require any  approximation of smooth functions by analytic ones. 

\giu

The main other novelty is in the new choice of the subsequence of  the sequence $(q_n)$ of denominators of the best rational approximations along which the cheap trick is applied.
Given the sequence $(q_n)$  of denominators of the best rational approximations 
of an irrational number, the notion of Diophantine bridges, a sufficiently long succession of $q_n$ where $q_{n+1}$ is bounded by a fixed power of $q_n$, was explicitly introduced in 
\cite{FKhanin} to settle the global smooth conjugacy problem of commuting circle diffeomorphisms. 

In \cite{3}, a pattern of Diophantine bridges that alternate with big power jumps in the sequence $(q_n)$ was proved to exist for every irrational $\a$. A subsequence $(q_{n_h})$ (called $(Q_h)$ in \cite{3}) of $(q_n)$ is then chosen using the endpoints of the Diophantine bridges and the $q_n$ with a big jump from $q_n$ to $q_{n+1}$. This subsequence is used to apply the cheap trick and prove rotations-reducibility in the real analytic context.   

The choice of the sequence  $(q_{n_h})$ is important to have a nice control on the Birkhoff sums and make sure that they are close to their mean (which allows to use the hypothesis on the fibered rotation number), and also to get sufficient gain from applying the cheap trick in terms of decreasing the magnitude of the non-abelian part of the cocycle. 

The choice of  the sequence  $(q_{n_h})$ in \cite{3} was adapted to the real analytic category, and in the current work it is crucial to choose quite differently the sequence  $(q_{n_h})$ to adapt to the smooth setting, that are needed not only to have a nice control on the Birkhoff sums, but also to have a good control of high norms of the error term.

\giu
Finally, to prove the convergence of the scheme it is crucial also to have improved estimates on the derivatives of the iterated cocycle that are adapted to the smooth setting (see Proposition \ref{p3}), compared to the previously used bounds such as the ones in  \cite{2}.

\giu
We now recall the idea behind the cheap trick and give an outline of the iterative conjugation scheme that is used to prove Theorem \ref{theorem1}.

\medskip 

Let $q_n$ be a denominator of a best rational approximation of $\a$. If the cocycle matrix $A(\cdot)=R_{\phi(\cdot)}+F(\cdot)$ is very close to the $SO(2,\Bbb{R})$ cocycle $R_{\phi(\cdot)}$, then $A^{(q_n)}(x):=R_{S_{q_n}\phi}(x)+\xi(x)$ is elliptic for every $x\in \mathbb{T}$.  In this case, we can find $B_1\in C^{\infty}(\Bbb{T},SL(2,\Bbb{R})),\phi_1\in C^{\infty}(\Bbb{T}, \Bbb{R})$ such that $$B_1(x)A^{(q_n)}(x)B_1^{-1}(x)=R_{\phi_1(x)},$$ and for $t\in\Bbb{N}$
\begin{equation} \label{intro_eq1}
\parallel B_1-Id.\parallel_t,\parallel\phi_1-\phi \parallel_t\leq C(t)\parallel(R_{2S_{q_n}\phi}-Id.)^{-1}\xi R_{-S_{q_n}\phi}\parallel_t,
\end{equation}

If in addition the Birkhoff sum $S_{q_n}\phi$ is uniformly close to its average, then this average is close to $q_n \rho(\a,A)$, where $\rho(\a,A)$ is the fibered rotation number of $(\a,A)$. Hence,  if $\norm{2q_n \rho(\a,A)}$ ($\norm{\cdot}$ being the distance to the closest integer) is not too small, then $B_1$ is close to the Identity matrix. Finally, the essential point in the cheap trick is to see that the conjugated cocycle $A_1(x):=B_1(x+q_n\a)A^{(q_n)}(x)B_1^{-1}(x)$ is of the order of $\norm{q_n\a}\sim \su{q_{n+1}}$ closer to a rotation valued cocycle then $A^{(q_n)}$ was. One can conjugate $A_1$ now and gain another factor 
$\su{q_{n+1}}$. Repeating the procedure $r_0+1$ times one obtains 
a conjugacy $\tilde{B}:=B_{r_0+1}\dots B_1$ such that 

\beqano
\tilde{B}(x+q_n\a)A^{(q_n)}(x)\tilde{B}^{-1}(x)=\tilde{F}(x)+R_{\phi_{h}(x)}
\eeqano
such that, if we suppose $\norm{2q_n \rho(\a,A)}\geq \frac{c}{n^2}$, it holds for $0\leq h\leq r_0+1$:
\begin{equation} \label{eq.RF44} 
\parallel\tilde{F}\parallel_t\leq   \frac{C(t+h)n^{2r_0(t+1)}}{q_{n+1}^h}\max_{\b\in\{0,1\}}\parallel\bar{F}\parallel_{t+h}^{\b}\left(\parallel\bar{F}\parallel_1\parallel\bar{\phi}\parallel_{t+h}\right)^{1-\b}.
\end{equation}
One then derives from \eqref{eq.RF44} similar estimates on how close  $\tilde{B}$ conjugates the original cocycle $(\a,A)$ to a rotation valued cocycle. 

An important point of the cheap trick is that we get inverse powers $q_{n+1}^{-h}$ in the control of the new nonlinearity when we compare to   $\parallel {F}\parallel_{t+h}$. However, there is no loss of derivatives in the scheme because we never solve a cohomological equation. 

\bigskip 

\noindent {\bf Plan of the paper.}    The main inductive conjugation step is stated in Proposition \ref{proposition1} and the 
    proof of Theorem \ref{theorem1}, based on  Propositon \ref{proposition1}, are  given in Section \ref{sect.main}.

 The outcome of the cheap trick is the content of  Propositon \ref{p2} of Section \ref{sec.cheap}, which is the main step of the proof.   The proof of Proposition \ref{proposition1} is given in Section \ref{sec5}. It goes through the application of the cheap trick to the iterated cocycle $A^{(q_n)}$, explained in Section \ref{sec.iterate}, and through the estimates for going back to the original cocycle, explained in Section \ref{sec.back}.  

The estimates on the iterated cocycle $(q_n\a,A^{(q_n)})$ are included in Section \ref{sec.iterated}. In Section \ref{sec.iterated1} we state and prove the estimates on the upper bounds on the non-abelian part of the iterated cocycle, and in Section \ref{sec.iterated2} we include the crucial estimates on the Birkhoff sums along the selected subsequence of $(q_n(\a))$. 

Finally, the Appendix A contains the statements and proofs of some basic analysis lemmas that are used throughout the paper.

\section{Notations and definitions}

In all the sequel $\a$ will be a fixed irrational number. 
If $f\in C^{\infty}(\Bbb{T},\Bbb{R})$ we denote the Birkhoff sums of $f$ above the circle rotation of angle $\a$ for $n\in\Bbb{N}$
\beqano
S_nf(x):=\sum_{h=0}^{n-1}f(x+h\a).
\eeqano
For $A\in C^{\infty}(\Bbb{T},SL(2,\Bbb{R})), n\in\Bbb{N}$, we denote
\beqano
A^{(n)}(x):=\prod_{j=n-1}^0 A(x+j\a).
\eeqano
For $f\in C^{\infty}(\Bbb{T},\Bbb{R}), k\in\Bbb{N}$, we denote by $D^k f$ the $k$-th derivative of $f$ and
\beqano
|f|_0:=\sup_{x\in\Bbb{T}}|f(x)|, \quad \parallel F\parallel_k:=\sum_{h=0}^k |D^hf|_0.
\eeqano

We use the notation 
\begin{displaymath}
  \norm{x} = \inf_{p \in \Z} | x - p|.
\end{displaymath}

We define $(q_n)$ to be the sequence of denominators of the best rational approximations of $\a$. Recall that $q_n$ satisfies $q_0=1$ and 
  \begin{equation} \label{best} \forall 1 \leq k < q_n,\qquad \norm{k\a} \geq \norm{q_{n-1}\a}. \end{equation}
 
Recall also that 
\begin{equation}{1\over q_{k+1}+q_k}<\norm{q_k \a}<{1\over q_{k+1}}. \label{best2} \end{equation}

\section{Main inductive conjugation step} \label{sect.main} 
In this section we state the main Proposition \ref{proposition1} containing one step of the reducibility scheme and we show how Theorem \ref{theorem1} follows directly from it. 
\begin{definition} \label{definition1}
    Let $\a$ be irrational. We define a subsequence of convergents $(q_{n_h})_{h\in\Bbb{N}}$ in the following way. We let $n_0=0$.
    Now, suppose that we have defined the subsequence up to $q_{n_h}$. If there exists $k$ such that $q_{n_h+1}^2\leq q_k<q_{n_h+1}^4$, we define ${n_{h+1}}:=k$. Otherwise, we take $n_{h+1}:=\max \{k\in\Bbb{N}:q_k\leq q_{n_h+1}^2\}$. For $h\in\Bbb{N}$, we use the notation $s_h:=\prod_{l=0}^{h}q_{n_l}$.
\end{definition}
\lem{lb} \label{lemma1}
For $h\in\Bbb{N}$, $q_{n_{h+1}}\leq q_{n_{h}+1}^4$. Moreover, for $h\in\Bbb{N}$: $s_h^6\leq q_{n_{h+1}}^{12}$.
\elem
\proof
It follows by the definition of $(q_{n_h})_{h\in\Bbb{N}}$. \qed


\pro{p1} \label{proposition1}
Let $\e>0$. There exists $\e_0>0$ such that the following holds. Let  $\a$ be irrational, $(q_n)_{n\in\Bbb{N}}$ be the convergents of $\a$, $r_0\geq 200$ and $A_0=R_{\phi_0}+F_0\in C^{\infty}(\Bbb{T}, SL(2,\Bbb{R})), \phi_0 \in \Bbb{T}, R_{\phi_0} \in SO(2,\Bbb{R})$ such that 
\beqano
\parallel F_0\parallel_{50r_0}<\e_0.
\eeqano
and 
$\rho:=\rho(\a,A_0)$ satisfies $\norm{2q_{n_h}\rho}\geq \frac{\e}{n_h^2}$ where $(n_h)$ is as in Definition 1. 

Then, there exist $B_{h}, F_h\in C^{\infty}(\Bbb{T},SL(2, \Bbb{R}))$, $\phi_h\in C^{\infty}(\Bbb{T}, SL(2,\Bbb{R}))$ such that for $h\in\Bbb{N}$:
\beq{e1}
B_h(x+\a)(R_{\phi_h(x)}+F_h(x))B_h^{-1}(x)=R_{\phi_{h+1}(x)}+F_{h+1}(x), 
\eeq
\beq{e2}
\parallel F_{h+1}\parallel_{1}\leq \frac{\e_0}{q_{n_{h+1}}^{\frac{r_0}{8}}},\quad \parallel F_{h+1}\parallel_{50r_0},\parallel\phi_{h+1}\parallel_{50r_0}\leq s_h^{6}
\eeq{}
where $s_h:=\prod_{l=0}^{h}q_{n_l}$ and, for $t\in\Bbb{N}$:
\beq{e3}
\parallel B_h-Id.\parallel_t, \parallel\phi_{h+1}-\phi_{h}\parallel_t, \parallel F_{h+1}\parallel_t\leq C(t)n_{h}^{4tr_0}q_{n_{h}}^5 \max_{\b\in\{0,1\}}\parallel F_h\parallel_{t}^{\b}(\parallel F_{h}\parallel_1q_{n_h}\parallel\phi_h\parallel_t)^{1-\b}.
\eeq
\epro{}
Now we show how Theorem \ref{t1} follows directly from Proposition \ref{p1}.

\giu
\nl
{\bf{Proof of Theorem \ref{t1}.}}
We want to show that the following limits exist in the $C^\infty$ category and that they satisfy all the requirements of Theorem \ref{theorem1} 
\beqano
B(x):=\lim_{h\rightarrow+\infty}B_{h}\dots B_{0}(x), \quad \phi(x):=\lim_{h\rightarrow+\infty}\phi_h(x).
\eeqano
 We first address the convergence in class $C^{r_0}$. Proposition \ref{proposition1} implies by convexity that 
\beq{e4}
\parallel F_h\parallel_{r_0}\leq C\parallel F_h\parallel_0^{1-\su{50}}\parallel F_h\parallel_{50r_0}^{\su{50}}<\frac{\sqrt{\e_0}}{q_{n_{h}}^6},
\eeq where in the last inequality we have used the estimate in (\ref{e2}) for $\parallel F_{h}\parallel_{1}, \parallel F_{h}\parallel_{50r_0}$ and the fact that $r_0\geq 200$. Next, by (\ref{e3}), (\ref{e4}):
\beqano
\parallel B-Id\parallel_{r_0},\parallel\phi-\rho\parallel_{r_0}&\leq& \prod_{h=0}^{\infty}(1+C(r_0)n_h^{4r_0^2}q_{n_h}^5\sup_{\b\in[0,1]}\parallel F_h\parallel_{r_0}^{\b}(\parallel F_{h}\parallel_0q_{n_h}\parallel\phi_h\parallel_{r_0})^{1-\b})-1\\
&\leq& \prod_{h=0}^{\infty}\left(1+\frac{C(r_0)n_h^{4r_0^2}\e_0}{q_{n_h}}\right)-1\leq \sqrt{\e_0}
\eeqano
if $\e_0$ is small enough. So, we have proved the convergence in norm $C^{r_0}$ as well as the bounds \eqref{eq.small}. The conjugacy equation \eqref{eq.conj} then follows from  \eqref{e1}. 

Now we use interpolation and the bounds in \eqref{e3} to prove the smooth convergence. By (\ref{e2}), (\ref{e3}):
\begin{align*}
\parallel\phi_{h+1}\parallel_t&\leq \parallel\phi_{h+1}-\phi_{h}\parallel_t+\parallel\phi_{h}\parallel_t\leq C(t)n_h^{4r_0^2}q_{n_h}^5(\parallel F_h\parallel_{r_0}+q_{n_h}\parallel F_{h}\parallel_0\parallel\phi_h\parallel_t)+\parallel\phi_h\parallel_t\\
&\leq C(t)n_h^{4r_0^2}q_{n_h}^5\parallel F_h\parallel_{r_0}+(1+\su{q_{n_h}})\parallel\phi_h\parallel_t.
\end{align*} 
So, iterating we get:
\beq{e5}
\parallel\phi_{h+1}\parallel_t\leq C(t)n_h^{4r_0^2}q_{n_h}^5\sum_{l=0}^h\parallel F_l\parallel_t.
\eeq{}

\giu
\nl
Let $t\in\Bbb{N}, t>r_0$. By (\ref{e3}) and (\ref{e5}), there exists $N\in\Bbb{N}$ such that for $h\geq N$:
\beqano
\parallel F_{h+1}\parallel_{20t}\leq n_h^{4r_0^2}q_{n_h}^5\sum_{l=0}^h\parallel F_l\parallel_{20t}\leq q_{n_{h+1}}^4 \sum_{l=0}^h\parallel F_l\parallel_{20t},
\eeqano
where in the last inequality we have used Lemma \ref{lb}. In particular there exists $N_1>N$ such that, for $h>N_1$:
\beqano
\parallel F_{h}\parallel_{20t}\leq n_h\left(\prod_{l=0}^{h-N}q_{n_{N+l}}^4\right)\sum_{l=0}^N\parallel F_l\parallel_{20t}\leq q_{n_{h+1}}^9\sum_{l=0}^N\parallel F_l\parallel_{20t}\leq q_{n_{h+1}}^{10}.
\eeqano
So, by convexity and the estimates of $\parallel F_h\parallel_{r_0}$ in (\ref{e4}) we get:
\beqano
\parallel F_h\parallel_t\leq C\parallel F_h\parallel_{r_0}^{1-\frac{t-r_0}{20t-r_0}}\parallel F_h\parallel_{20t}^{\frac{t-r_0}{20t-r_0}}\leq \frac{1}{q_{n_{h}}^{5+\su{3}}}.
\eeqano
In particular, for any $t\in\Bbb{N}$ there exists $\bar{N}\in\Bbb{N}$ such that for $h>\bar{N}$:
\beqano
\parallel B_h-Id.\parallel_t,\parallel\phi_{h+1}-\phi_h\parallel_t\leq \frac{C(t)n_h^{4tr_0}}{q_{n_h}^{\su{3}}},
\eeqano
finishing the proof of smooth convergence of the scheme. \qed

\section{The cheap trick} \label{sec.cheap} 
The aim of this section is to prove the following.
\pro{p2}
Let ${\bar{\phi}}  \in C^{\infty}(\Bbb{T}, \Bbb{R})$ and $\bar{A}=R_{\bar{\phi}}+\bar{F}\in C^{\infty}(\Bbb{T}, SL(2,\Bbb{R})),r_0\in\Bbb{N}, n\in \Bbb{N}$ be such that 
\begin{equation} \label{eq.RF} 
\parallel(R_{2\bar{\phi}}-Id.)^{-1}\parallel_0\leq Cn^2,\quad  \parallel\bar{F}\parallel_0<\su{q_n}, \parallel\bar{F}\parallel_{r_0},\parallel\bar{\phi}\parallel_{r_0}<1.
\end{equation} 
Then, there exist $\tilde{B},\tilde{F}\in C^{\infty}(\Bbb{T},SL(2,\Bbb{R})), \tilde{\phi}\in C^{\infty}(\Bbb{T}, \Bbb{R})$ with for any $t\in \Bbb{N}$
\begin{equation} \label{eq.RF2} 
\parallel\tilde{B}-Id.\parallel_t,\parallel\tilde{\phi}-\bar{\phi}\parallel_t\leq C(t)n^{2(t+1)}\max_{\b\in\{0,1\}}\parallel\bar{F}\parallel_t^{\b}\left(\parallel\bar{F}\parallel_1\parallel\bar{\phi}\parallel_t\right)^{1-\b}
\end{equation} 
and
\begin{equation} \label{eq.RF3} 
\tilde{B}(x+q_n\a)\bar{A}(x)\tilde{B}^{-1}(x)=R_{\tilde{\phi}}+\tilde{F},
\end{equation}
and, for $0\leq h\leq r_0+1$:
\begin{equation} \label{eq.RF4} 
\parallel\tilde{F}\parallel_t\leq   \frac{C(t+h)n^{2r_0(t+1)}}{q_{n+1}^h}\max_{\b\in\{0,1\}}\parallel\bar{F}\parallel_{t+h}^{\b}\left(\parallel\bar{F}\parallel_1\parallel\bar{\phi}\parallel_{t+h}\right)^{1-\b}.
\end{equation}
\epro
In order to prove Proposition \ref{p2} we prove at first some simple lemmas. 

\giu

\lem{l5}
Let $D>0$. There exists $\e>0$ such that, if $\bar{A}=R_{\bar{\phi}}+\bar{F}\in C^{\infty}(\Bbb{T},SL(2,\Bbb{R}))$ with $\parallel A\parallel_0<D,\parallel\bar{F}\parallel<\e\min\{1, \parallel R_{2\bar{\phi}}-Id.\parallel_0\}$, then there exists $\tilde{B}_1\in C^{\infty}(\Bbb{T}, SL(2,\Bbb{R})), \tilde{\phi}_1\in C^{\infty}(\Bbb{T},\Bbb{R})$ with:
\beq{e10}
\parallel\tilde{B}_1-Id.\parallel_t,\parallel\bar{\phi}-\tilde{\phi}_1\parallel_t\leq C(t)\parallel(R_{2\bar{\phi}}-Id.)^{-1}\bar{F}R_{-\bar{\phi}}\parallel_t
\eeq
such that:
\beq{e11}
\tilde{B}_1(x)\bar{A}(x)\tilde{B}_1^{-1}(x)=R_{\tilde{\phi}_1(x)}.
\eeq{}
\elem{}
\proof
We will sketch the proof following \cite{3}. Let $Y:=\bar{F}R_{-\bar{\phi}}$, $G_1:=\log(1+Y), \theta_1=\bar{\phi}$, so that $\bar{A}=e^{G_1} R_{\theta_1}$.
Let:
\beqano
G_1&=&
\begin{pmatrix}
x_1 & y_1-2\pi z_1 \\
y_1+2\pi z_1 & -x_1
\end{pmatrix},\\
\bar{G}_1&:=&
\begin{pmatrix}
x_1 & y_1 \\
y_1 & -x_1
\end{pmatrix}, \\
\begin{pmatrix}
\tilde{x}_1 \\
\tilde{y}_1
\end{pmatrix}
&:=& (R_{2\theta_1}-Id.)^{-1}
\begin{pmatrix}
x_1 \\
y_1
\end{pmatrix},\\
v_1&:=&
\begin{pmatrix}
\tilde{x}_1 & \tilde{y}_1 \\
\tilde{y}_1 & -\tilde{x}_1
\end{pmatrix},
\quad \theta_2:=\theta_1+z_1.
\eeqano
Then, we want to show that:
\beqano
e^{v_1}e^{G_1}R_{\theta_1}e^{-v_1}=e^{G_2}R_{\theta_2},
\eeqano
with:
\beq{e7}
\parallel e^{v_1}-Id.\parallel_t\leq \parallel(R_{2\theta_1}-Id.)^{-1}G_1R_{-\theta_1}\parallel_t,\quad \parallel\theta_2-\theta_1\parallel_t\leq \parallel G_1\parallel_t,
\eeq{}
\beq{e8}
\parallel G_2\parallel_t\leq C(t)\parallel(R_{2\theta_1}-Id.)^{-1}G_1^2R_{-\theta_1}\parallel_t.
\eeq{}
Note that (\ref{e7}) follows directly from the definition of $v_1, \theta_1$. So, we just prove (\ref{e8}):
\beqano
e^{v_1}e^{G_1}R_{\theta_1}e^{-v_1}=(Id.+v_1)(Id.+G_1)R_{\theta_1}(Id.-v_1)+O(v_1^2).
\eeqano
Then,  since by definition of $v_1$ we have
\beq{e6}
v_1R_{\theta_1}-R_{\theta_1}v_1+\bar{G}_1R_{\theta_1}=0,
\eeq{}
we get \beqano
(Id.+v_1)(Id.+G_1)R_{\theta_1}(Id.-v_1)&=&(Id.+G_1-\bar{G}_1)R_{\theta_1}+v_1R_{\theta_1}\\
&-&R_{\theta_2}v_1+\bar{G}_1R_{\theta_1}+O(v_1^2)\\
&=&(Id.+G_1-\bar{G}_1)R_{\theta_1}+O(v_1^2)\\
&=&R_{\theta_2}+O(v_1^2),
\eeqano
where the last inequality follows from \beq{e9}
Id.+G_1-\bar{G}_1-R_{\theta_2-\theta_1}=O(G_1^2).
\eeq{}
Next, (\ref{e8}) follows by the first inequality ($\parallel e^{v_1}-Id.\parallel_t\leq \parallel(R_{2\theta_1}-Id.)^{-1}G_1R_{-\theta_1}\parallel_t$) and by Lemma \ref{v1} of the Appendix.
Finally, the Lemma follows by iterating the scheme, with $\tilde{B}_1:=\lim e^{v_n}\dots e^{v_1}$, $\tilde{\phi}_1=\lim \theta_n$.
\qed
\lem{l7}
Let $\phi\in C^{\infty}(\Bbb{T},\Bbb{R})$ such that $\parallel R_{2\phi}-Id.\parallel_0\geq\frac{C}{n^2}$. For $t\in\Bbb{N}$:
\beqano
|D^{t}(R_{2\phi}-Id.)^{-1}|_0\leq C(t)n^{2(t+1)}\parallel\phi\parallel_t
\eeqano
\elem{}\proof
It follows by Lemma \ref{l6} of the Appendix by applying Hadamard's inequality to each term of the homogeneous polynomials. \qed

\nl
{\bf{Proof of Proposition \ref{p2}.}}
By \eqref{eq.RF}, we have that $\parallel A\parallel_0<D,\parallel\bar{F}\parallel<\e\min\{1, \parallel R_{2\bar{\phi}}-Id.\parallel_0\}$. ï¿½Hence, Lemma \ref{l5} gives $\tilde{B}_1, \tilde{\phi}_1$ such that (\ref{e10}) and (\ref{e11}) hold.
Then, by \eqref{e10} of Lemma \ref{l5}, and Lemma \ref{l7}  and Leibnitz formula, we get for $h\in\Bbb{N}$
\beqano
|D^h(\tilde{B}_1-Id)|_0, \leq C(h)n^{2(h+1)}\sum_{h_1+h_2+h_3=h}|D^{h_1}\bar{F}|_0\parallel\bar{\phi}\parallel_{h_2}\parallel\bar{\phi}\parallel_{h_3}.
\eeqano
Hence,  by interpolation inequalities (see Lemma \ref{v5} in the Appendix), we get  for $t\in\Bbb{N}$:\beq{w}
|D^t(\tilde{B}_1-Id.)|_0\leq C(t)n^{2(t+1)}\max_{\b\in\{0,1\}}\parallel\bar{F}\parallel_t^{\b}\left(\parallel\bar{F}\parallel_0\parallel\bar{\phi}\parallel_t\right)^{1-\b} .\eeq
In the same way, \eqref{e10} implies 
\beqano\parallel\bar{\phi}-\phi_1\parallel_t\leq C(t)n^{2(t+1)}\max_{\b\in\{0,1\}}\parallel\bar{F}\parallel_{t}^{\b}\left(\parallel\bar{F}\parallel_0\parallel\bar{\phi}\parallel_{t}\right)^{1-\b}.\eeqano
Now, let
\beqano
F_1(x)=\tilde{B}_1(x+q_n\a)\bar{A}(x)\tilde{B}_1^{-1}(x)-R_{\phi_1(x)}=(\tilde{B}_1(x+q_n\a)-\tilde{B}_1(x))\bar{A}(x)\tilde{B}_1^{-1}(x).
\eeqano
We want to show that, for $l=0,1$:
\beq{w1}
\parallel F_1\parallel_t\leq   \frac{C(t)n^{2(t+1+l)}}{q_{n+1}^l}\max_{\b\in\{0,1\}}\parallel\bar{F}\parallel_{t+l}^{\b}\left(\parallel\bar{F}\parallel_0\parallel\bar{\phi}\parallel_{t+l}\right)^{1-\b}
\eeq{}
We prove at first (\ref{w1}) for $l=1$.

\giu
Note that:
$|D^{t}(\tilde{B}_1(x+q_n\a)-\tilde{B}_1(x))|_0\leq \frac{1}{q_{n+1}}|D^{t+1}(\tilde{B}_1-Id.)|_0$. 
Then:
\beqano
|D^t F_1|_0 &\leq& \sum_{t_1+t_2+t_3=t}|D^{t_1}(\tilde{B}_1(x+q_n\a)-\tilde{B}_1(x))|_0|D^{t_2}\bar{A}|_0|D^{t_3}\tilde{B}_1^{-1}|_0\\
&\leq& \su{q_{n+1}}\sum_{t_1+t_2+t_3=t}|D^{t_1+1}(\tilde{B}_1-Id.)|_0|D^{t_2}\bar{A}|_0|D^{t_3}\tilde{B}_1^{-1}|_0\\
&\leq&  \frac{C(t)n^{2(t+2)}}{q_{n+1}^h}\max_{\b\in\{0,1\}}\parallel\bar{F}\parallel_{t+1}^{\b}\left(\parallel\bar{F}\parallel_0\parallel\bar{\phi}\parallel_{t+1}\right)^{1-\b},
\eeqano
with the last inequality that follows from (\ref{w}).

Now we prove (\ref{w1}) for $l=0$:
\beqano
|D^t F_1|_0 &\leq& \sum_{t_1+t_2+t_3=t,t_1>0}|D^{t_1}(\tilde{B}_1(x+q_n\a)-\tilde{B}_1(x))|_0|D^{t_2}\bar{A}|_0|D^{t_3}\tilde{B}_1^{-1}|_0\\
&+& \sum_{t_2+t_3=t}|(\tilde{B}_1(x+q_n\a)-\tilde{B}_1(x))|_0|D^{t_2}\bar{A}|_0|D^{t_3}\tilde{B}_1^{-1}|_0.\\
\eeqano
For the first term we get:
\beqano
\sum_{t_1+t_2+t_3=t,t_1>0}|D^{t_1}(\tilde{B}_1(x+q_n\a)-\tilde{B}_1(x))|_0|D^{t_2}\bar{A}|_0|D^{t_3}\tilde{B}_1^{-1}|_0
\leq
\eeqano
\beqano
2\sum_{t_1+t_2+t_3=t}|D^{t_1}\tilde{B}_1|_0|D^{t_2}\bar{A}|_0|D^{t_3}\tilde{B}_1^{-1}|_0\leq C(t)n^{2(t+2)}\max_{\b\in\{0,1\}}\parallel\bar{F}\parallel_{t}^{\b}\left(\parallel\bar{F}\parallel_0\parallel\bar{\phi}\parallel_{t}\right)^{1-\b}.
\eeqano
For the second term we get:
\beqano
\sum_{t_2+t_3=t}|(\tilde{B}_1(x+q_n\a)-\tilde{B}_1(x))|_0|D^{t_2}\bar{A}|_0|D^{t_3}\tilde{B}_1^{-1}|_0\leq
\eeqano
\beqano
\su{q_{n+1}}\sum_{t_2+t_3=t}|D(\tilde{B}_1-Id.)|_0|D^{t_2}\bar{A}|_0|D^{t_3}\tilde{B}_1^{-1}|_0\leq C(t)n^{2(t+2)}\max_{\b\in\{0,1\}}\parallel\bar{F}\parallel_{t}^{\b}\left(\parallel\bar{F}\parallel_1\parallel\bar{\phi}\parallel_{t}\right)^{1-\b}.
\eeqano

\giu
\nl
Now, suppose that we iterated the scheme $m$ times, with $m\leq r_0$. For $k\leq m$ let: 
\beqano
\tilde{\bar{B}}_k&=&\tilde{B}_{k}\dots\tilde{B}_1,\\
R_{\tilde{\phi}_k(x)}&=&\tilde{\bar{B}}_k(x)\bar{A}(x)\tilde{\bar{B}}_k^{-1}(x),\\
\tilde{F}_k(x)&=&(\tilde{\bar{B}}_k(x+q_n\a)-\tilde{\bar{B}}_k(x))\bar{A}(x)\tilde{\bar{B}}_k^{-1}(x).
\eeqano
Then, in the same way as above, for $l=0,1$ and $k\leq m$:
\beq{ii}
\parallel \tilde{F}_{k} \parallel_t, \parallel \tilde{\phi}_k-\tilde{\phi}_{k-1}\parallel_t\leq   \frac{C(t)n^{2(t+1+l)}}{q_{n+1}^l}\max_{\b\in\{0,1\}}\parallel \tilde{F}_{k-1}\parallel_{t+l}^{\b}\left(\parallel \tilde{F}_{k-1}\parallel_1\parallel\tilde{\phi}_{k-1}\parallel_{t+l}\right)^{1-\b},
\eeq{}
In particular, from $(\ref{ii})$ and (\ref{eq.RF}) we get:
\beqano
\parallel \tilde{F}_k\parallel_{r_0}, \parallel \tilde{\phi}_k\parallel_{r_0}\leq C(r_0)n^{4(r_0+1)k}, 
\eeqano
\beqano
\parallel \tilde{F}_k\parallel_{1},\parallel \tilde{\phi}_k-\tilde{\phi}_{k-1}\parallel_{1}&\leq& \frac{n^{4r_0}}{q_{n+1}}\max_{\b\in\{0,1\}}\parallel \tilde{F}_{k-1}\parallel_{r_0}^{\b}\left(\parallel \tilde{F}_{k-1}\parallel_1\parallel\tilde{\phi}_{k-1}\parallel_{r_0}\right)^{1-\b}\\
&\leq& \frac{C(r_0)n^{4(r_0+1)(k+1)}}{q_{n+1}}<1.
\eeqano
So, because $\tilde{B}_{m}(x+\a)A(x)\tilde{B}_m(x)^{-1}=R_{\tilde{\phi}_m(x)}+\tilde{F}_m(x)$ satisfies:
\beqano
\parallel R_{\tilde{\phi}_m}+\tilde{F}_m\parallel_0<2<D, \parallel\tilde{F}_m\parallel_0\leq \e\min\{1, \parallel R_{2\tilde{\phi}_m}-Id.\parallel_0\},
\eeqano
we can iterate the scheme one more time. Finally, Proposition \ref{p2} follows with:
\beqano
\tilde{B}:=\tilde{B}_{r_0+1}\dots\tilde{B}_1, \tilde{\phi}:=\tilde{\phi}_{r_0+1}, \tilde{F}:=\tilde{F}_{r_0+1}.
\eeqano
\qed
\section{Estimates for the iterated cocycle} \label{sec.iterated}Ê
Now we start by giving some bounds for the derivatives of the iterates of the cocycle.

\giu
In the following proposition we will show that if $A=R_{\phi}+F$ is quite close to a rotation valued cocycle ($F$ small), then also $A^{(q_n)}$ remains close to the rotation valued cocycle $R_{S_{q_n}\phi}$. In Section \ref{sec.iterated2}, we will see that $S_{q_n}\phi$ is close to its average (that is close to $q_n\rho$ with $\rho=\rho(\a,A)$).
\subsection{Estimates of the derivatives of the iterated cocycle }  \label{sec.iterated1}Ê
The goal of this section is to prove the following 
\pro{p3}
Let $A=R_{\phi}+F\in C^{\infty}(\Bbb{T}, SL(2,\Bbb{R}))$ be such that:
\beq{}
\parallel F\parallel_{0}<\su{q_n^6}, \quad \parallel\phi-\hat{\phi}(0)\parallel_0<1.
\eeq
Let $\xi$ be defined such that $A^{(q_n)}=R_{S_{q_n}\phi}+\xi$. Then, for $t\in\Bbb{N}$:
\beq{}
\parallel\xi\parallel_t\leq C(t)q_n^5\left(\parallel F\parallel_t+\parallel F\parallel_0\max_{\b\in\{0,1\}}\parallel F\parallel_t^{\b}\parallel\phi\parallel_t^{1-\b}\right).
\eeq
\epro{}

\nl

\nl
{\bf{Proof of Proposition 3. }}
Write $\xi:=A^{(q_n)}-R_{S_{q_n}\phi}$ as:
\beqano
\xi(x)=\sum_{h=1}^{q_n}\sum_{q_n-1\geq i_1>i_2>...>i_h\geq 0}\left(\prod_{j=q_n-1}^{i_1+1}R_{\phi(x+j\a)}\right)F(x+i_1\a)\dots F(x+i_h\a)\left(\prod_{j=i_{h}-1}^{0}R_{\phi(x+j\a)}\right),
\eeqano
Then:
\begin{multline*} \xi(x)=\\
\sum_{h=1}^{q_n}\sum_{q_n-1\geq i_1>i_2>...>i_h\geq 0}R_{S_{(q_n-i_1-1)}\phi(x+(i_1+1)\a)}F(x+i_1\a)R_{S_{(i_1-i_2-1)}\phi(x+(i_2+1)\a)}F(x+i_2\a)\dots R_{S_{(i_h-1)}\phi(x)}.
\end{multline*}
So, let $t\in\Bbb{N}$. For each $1\leq h\leq q_n$ we have the sum of $\binom{q_n}{h}$ terms on the form:
\beqano
R_{S_{j_1}\phi(x+(i_1+1)\a)}F(x+i_1\a)\dots F(x+i_{h}\a)R_{S_{j_{h+1}}\phi(x)},
\eeqano
for some $0\leq j_l\leq q_n$ and for $1\leq l\leq h+1$ (and with the same convention $R_{S_{0}\phi}:=Id$). In particular, by Leibnitz formula, the $t$-h derivative of each of these terms is the sum of at most $(2h+1)^t$ terms of the form:
\beqano
D^{t_1}R_{S_{j_1}\phi(x+(i_1+1)\a)}D^{t_2}F(x+i_1\a)\dots D^{t_{2h}}F(x+i_{h}\a)D^{t_{2h+1}}R_{S_{j_{h+1}}\phi(x+(i_h-1)\a)},
\eeqano
with $t_1+\dots t_{2h+1}=t$. If $t_{2l}\not=0$, by the interpolation Lemma \ref{v5} of the Appendix, we get
\beqano
|D^{t_{2l}}F|_0\leq C\parallel F\parallel_{0}^{1-\frac{t_{2l}}{t}}\parallel F\parallel_t^{\frac{t_{2l}}{t}}.
\eeqano
In the same way, if $t_{2l+1}\not=0$, then by Lemma \ref{v6}:
\beqano
|D^{t_{2l+1}}R_{S_{j_{2l+1}}\phi}|_0\leq C(t)q_n\parallel\phi\parallel_{t_{2l+1}}\leq C(t)q_n\parallel\phi\parallel_{0}^{1-\frac{t_{2l+1}}{t}}\parallel\phi\parallel_t^{\frac{t_{2l+1}}{t}}.
\eeqano
Moreover, by the assumptions we have $\parallel\phi-\hat{\phi}(0)\parallel_0<1$.
So, it follows that:
\beqano
|D^t F|_{0}\leq C(t)q_n^3\parallel\phi\parallel_0\parallel F\parallel_t+\sum_{h=2}^{q_n}\binom{q_n}{h}(2h+1)^t q_{n}^{h+1} C(t)^{2h+1}\parallel F\parallel_{0}^{h-1}\sup_{\a\in[0,1]}\parallel F\parallel_{t}^{\a}\parallel\phi\parallel_{t}^{1-\a}
\eeqano
\beqano
\leq C(t)q_n^3\parallel F\parallel_t+\parallel F\parallel_0\sup_{\a\in[0,1]}\parallel F\parallel_{t}^{\a}\parallel\phi\parallel_{t}^{1-\a}\sum_{h=2}^{q_n} \frac{q_n^h}{h!}(2h+1)^t q_{n}^{h+1} C(t)^{2h+1}\parallel F\parallel_{0}^{h-2}
\eeqano
\beqano
\leq C(t)q_n^3\parallel F\parallel_t+\parallel F\parallel_0\sup_{\a\in[0,1]}\parallel F\parallel_{t}^{\a}\parallel\phi\parallel_{t}^{1-\a}\sum_{h=2}^{q_n} \frac{q_n^h}{h!}(2h+1)^t q_{n}^{h+1} C(t)^{2h+1}\parallel F\parallel_{0}^{h-2}.
\eeqano
Now, because of $\parallel F\parallel_0<\su{q_n^{6}}$, for $h\geq 2$ we get:
\beqano
q_n^{2h+1}\parallel F\parallel_{0}^{h-2}\leq \frac{q_n^5}{q_n^{h}}.
\eeqano
So:
\beqano
\sum_{h=2}^{q_n} \frac{q_n^h}{h!}(2h+1)^t q_{n}^{h+1} C(t)^{2h+1}\parallel F\parallel_{0}^{h-2}\leq q_n^5\parallel|\phi\parallel_{t}\sum_{h=0}^{+\infty}\left(\frac{C(t)}{q_n}\right)^h\su{h!}\leq C(t) q_n^5.
\eeqano
\qed
\subsection{Estimates of the Birkhoff sums} \label{sec.iterated2}
The goal of this section is to prove the following 

\medskip
\pro{ll8}
Let $\xi:=A^{(q_n)}-R_{S_{q_n}\phi}$ and suppose that $\norm{2q_n\rho}\geq\frac{\e}{n^2}, \parallel\phi\parallel_7\leq 1, \parallel\xi\parallel_0<\su{q_n}$. If $q_{n+1}>q_n^2$ or there exists $k$ such that $q_k^2<q_n<q_k^4$ then, for $t\geq 0$: 
\beqano
|D^{t}(R_{2S_{q_n}\phi}-Id.)^{-1}|_0\leq C(t)n^2\parallel S_{q_n}\phi\parallel_t.
\eeqano
\epro

We first show in Lemmas \ref{l10} and \ref{l11} that, under certain conditions on $q_n$, the Birkhoff sum $S_{q_n}\phi$ is close to its average. In the proofs, we will need Fourier truncation and rest operators that we now define. 
\begin{definition}
Let $f\in C^{\infty}(\Bbb{T}, \Bbb{R}), a>0$. 
\beqano
T_a(f):=\sum_{|l|\leq a}\hat{f}(l)e^{2\pi i lx}, \quad R_a(f):=\sum_{|l|> a}\hat{f}(l)e^{2\pi i lx}.
\eeqano
\end{definition}

\lem{l10}
Suppose that there exists $q_k$ such tha $q_k^2<q_n<q_k^4$. Then, for $t\in\Bbb{N}$:
\beqano
|D^t (S_{q_n}\phi-q_n\hat{\phi}(0)|_0\leq C(t)\min\{q_n\parallel\phi-\hat{\phi}(0)\parallel_t, \su{q_n^{\su{4}}}\parallel\phi-\hat{\phi}(0)\parallel_{t+7}\}.
\eeqano
\elem{}
\proof
Let $t\in\Bbb{N}$. 
The inequality 
\beqano
|D^t (S_{q_n}\phi-q_n\hat{\phi}(0)|_0\leq q_n C(t)\parallel\phi-\hat{\phi}(0)\parallel_t
\eeqano
is trivial.
Then, for the second inequality:
\beqano
|D^t(S_{q_n}\phi-q_n\hat{\phi}(0))|_0\leq |D^t T_{q_n^{\su{4}}}(S_{q_n}\phi-q_n\hat{\phi}(0))|_0+|D^t R_{q_n^{\su{4}}}(S_{q_n}\phi-q_n\hat{\phi}(0))|_0
\eeqano
So: 
\beqano
|D^t T_{q_n^{\su{4}}}(S_{q_n}\phi-q_n\hat{\phi}(0))|_0\leq \sum_{1\leq |l|\leq q_n^{\su{4}}}|2\pi l|^t\left|\frac{e^{2\pi i q_nl\a}-1}{e^{2\pi il\a}-1}\right||l\hat{\phi}(l)|
\eeqano
\beqano
\leq \frac{q_k}{q_{n+1}} \sum_{1\leq |l|\leq q_n^{\su{4}}}|2\pi l|^{t+1}|\hat{\phi}(l)|,
\eeqano
where in the last inequality we have used the fact that for $|l|\leq q_n^{\su{4}}$, $|e^{2\pi i q_nl\a}-1|\leq \frac{2\pi|l|}{q_{n+1}}$ and $|e^{2\pi il\a}-1|\geq \frac{1}{q_k}$ (that follows by $|l|\leq q_n^{\su{4}}<q_k$). Finally:
\beqano
\frac{q_k}{q_{n+1}}\sum_{1\leq |l|\leq q_n^{\su{4}}}|2\pi l|^{t+1}|\hat{\phi}(l)|\leq \frac{Cq_n^{\su{2}}}{q_{n+1}}\parallel\phi-\hat{\phi}(0)\parallel_{t+3},
\eeqano
with the last inequality that follows by $q_k\leq q_n^{\su{2}}$ and $|\hat{\phi}(l)|\leq \frac{\parallel\phi-\hat{\phi}(0)\parallel_{t+3}}{|2\pi l|^{t+3}}$. Moreover:
\beqano
|D^t R_{q_n^{\su{4}}}(S_{q_n}\phi-q_n\hat{\phi}(0))|_0\leq q_n|D^t R_{q_n^{\su{4}}}(\phi-\hat{\phi}(0))|_0\leq \frac{q_n}{q_n^{\frac{5}{4}}}\parallel\phi-\hat{\phi}(0)\parallel_{t+7},
\eeqano
with the last inequality that follows by Lemma \ref{l9}. Then:
\beqano
|D^t(S_{q_n}\phi-q_n\hat{\phi}(0))|_0\leq |D^t T_{q_n^{\su{4}}}(S_{q_n}\phi-q_n\hat{\phi}(0))|_0+|D^t R_{q_n^{\su{4}}}(S_{q_n}\phi-q_n\hat{\phi}(0))|_0
\eeqano
\beqano
\leq\frac{Cq_n^{\su{2}}}{q_{n+1}}\parallel\phi-\hat{\phi}(0)\parallel_{t+3}+\frac{q_n}{q_n^{\frac{5}{4}}}\parallel\phi-\hat{\phi}(0)\parallel_{t+7}\leq \frac{C}{q_n^{\su{4}}}\parallel\phi-\hat{\phi}(0)\parallel_{t+7}.
\eeqano
\qed
\lem{l11}
If $q_{n+1}>q_n^2$, then for $t\in\Bbb{N}$:
\beqano
|D^t( S_{q_n}\phi-q_n\hat{\phi}(0))|_0\leq \frac{C}{q_n}\parallel\phi-\hat{\phi}(0)\parallel_{t+4}.
\eeqano
\elem{}
\proof
Let $t\in\Bbb{N}$:
\beqano
|D^{t}(S_{q_n}\phi-\hat{\phi}(0))|_0\leq |D^t T_{q_n-1}(S_{q_n}\phi-q_n\hat{\phi}(0))|_0+|D^t R_{q_n-1}(S_{q_n}\phi-q_n\hat{\phi}(0))|_0
\eeqano
Then:
\beqano
|D^t T_{q_n-1}(S_{q_n}\phi-q_n\hat{\phi}(0))|_0\leq \sum_{1\leq |l|< q_n}|2\pi l|^t\left|\frac{e^{2\pi i q_nl\a}-1}{e^{2\pi il\a}-1}\right| |\hat{\phi}(l)|
\eeqano
\beqano
\leq \frac{Cq_n}{q_{n+1}}\sum_{1\leq |l|< q_n}|2\pi l|^{t+1}\frac{|D^{t+3}\phi|_0}{|l|^{t+3}}\leq \frac{C}{q_{n}}\parallel\phi-\hat{\phi}(0)\parallel_{t+3}.
\eeqano
Now we estimate the second term:
\beqano
|D^t R_{q_n-1}(S_{q_n}\phi-q_n\hat{\phi}(0))|_0\leq q_n|D^{t}R_{q_n-1}(\phi-\hat{\phi}(0))|_0
\eeqano
\beqano
\leq \frac{Cq_n}{q_n^2}\parallel\phi-\hat{\phi}(0)\parallel_{t+4},
\eeqano
with the last inequality that follows by Lemma \ref{l9}. \qed

\noindent {\bf Proof of Proposition \ref{ll8}. } 
Let $\xi:=A^{(q_n)}-R_{S_{q_n}\phi}$. By the definition of the fibered rotation number, we have that  
\beq{e21}
|\rho(R_{S_{q_n}\phi}+\xi)-q_n\hat{\phi}(0)|\leq |R_{S_{q_n}\phi}-R_{q_n\hat{\phi}(0)}|_0+|\xi|_0. 
\eeq{}
Moreover, by Lemma \ref{l10} and Lemma \ref{l11}:
\beqano
|R_{S_{q_n}\phi}-R_{q_n\hat{\phi}(0)}|_0\leq \frac{C\parallel\phi\parallel_7}{q_n^{\su{4}}}\leq  \frac{C}{q_n^{\su{4}}},
\eeqano
Then, the proposition follows by Lemma \ref{l7}, (\ref{e21}) and the assumptions $\parallel\xi\parallel_0<\su{q_n}$ and $\norm{2q_n\rho}\geq\frac{\e}{n^2}$. \qed
\section{Proof of Proposition \ref{p1}} \label{sec5}
\nl

\subsection{Applying the cheap trick to the iterated cocycle}Ê\label{sec.iterate}ÊThe estimates in (\ref{e2}) for $h=0$ are trivially satisfied by the hypothesis in Theorem \ref{t1}. 
Now, suppose that we are at the $h$-step so that we have the cocycle $A_h=R_{\phi_h}+F_h$ such that the estimates in (\ref{e2}) hold, that is \beq{e16}
\parallel F_{h}\parallel_{1}\leq \frac{\e_0}{q_{n_{h}}^{\frac{r_0}{8}}},\quad \parallel F_{h}\parallel_{50r_0},\parallel\phi_{h}\parallel_{50r_0}\leq s_h^6,
\eeq{}
with $s_h=\prod_{l=0}^{h}q_{n_l}$. We want to find $B_{h}, F_{h+1}\in C^{\infty}(\Bbb{T}, SL(2,\Bbb{R})),\phi_{h+1}\in C^{\infty}(\Bbb{T},\Bbb{R})$ such that (\ref{e1}), (\ref{e2}), (\ref{e3}) hold. We have 
\pro{p4}
There exist $B_{h},\tilde{F}\in C^{\infty}(\Bbb{T},SL(2,\Bbb{R})), \tilde{\phi}\in C^{\infty}(\Bbb{T}, \Bbb{R})$ with:
\beq{lp1}
\parallel B_h-Id.\parallel_t,\parallel\tilde{\phi}-\bar{\phi}\parallel_t\leq C(t)n^{2r_0(t+1)}q_{n_h}^5\max_{\b\in\{0,1\}}\parallel F_h\parallel_t^{\b}\left(\parallel F_h\parallel_1\parallel\bar{\phi}\parallel_t\right)^{1-\b}
\eeq{}
such that:
\beq{lp2}
B_h(x+q_{n_h}\a)A^{(q_{n_h})}(x)B_h^{-1}(x)=R_{\tilde{\phi}}+\tilde{F},
\eeq{}
and, for $0\leq l\leq r_0+1$:
\beq{lp3}
\parallel\tilde{F}\parallel_t\leq   \frac{C(t+l)q_{n_h}^5n_h^{2r_0(t+2)}}{q_{n_h+1}^l}\max_{\b\in\{0,1\}}\parallel F_h\parallel_{t+l}^{\b}\left(\parallel F_h\parallel_1 q_{n_h}\parallel\phi_h\parallel_{t+l}\right)^{1-\b}.
\eeq{}
In particular:
\beq{ee18}
\parallel\tilde{F}\parallel_{1}\leq \frac{Cq_{n_h}^5 n_h^{4r_0}}{q_{n_{h+1}}^{r_0}}\max_{\b\in\{0,1\}}\parallel F_h\parallel_{r_0+1}^{\b}\left(\parallel F_h\parallel_1q_{n_h}\parallel\phi_h\parallel_{r_0+1}\right)^{1-\b}\leq\frac{Cq_{n_h}^5 n_h^{4r_0}}{q_{n_{h+1}}^{r_0}},
\eeq{}
and
\beq{lp4}
\parallel\tilde{F}\parallel_{50r_0},\parallel\tilde{\phi}-\bar{\phi}\parallel_{50r_0}\leq C(r_0)q_{n_h}^5 n_h^{2r_0(50r_0+2)}\max_{\b\in\{0,1\}}\parallel F_h\parallel_{50r_0}^{\b}\left(\parallel F_h\parallel_1 q_{n_h}\parallel\phi_h\parallel_{50r_0}\right)^{1-\b}.
\eeq{}

Moreover, 
\beq{coro1}
\parallel\tilde{F}\parallel_{50r_0},\parallel\tilde{\phi}-\bar{\phi}\parallel_{50r_0}\leq C(r_0)q_{n_h}^{5} s_h^6 n_h^{2r_0(50r_0+2)}.
\eeq
\epro

\proof
By Lemma \ref{l10}, Lemma \ref{l11} we get:
\beq{e17}
\parallel S_{q_{n_h}}\phi_h-q_{n_h}\hat{\phi_h}(0)\parallel_{0}\leq \frac{C}{q_{n_h}^{\su{4}}}\parallel\phi_h-\hat{\phi_h}(0)\parallel_7\leq\su{q_{n_h}^{\su{4}}}.
\eeq{}
Then, by the arithmetic condition $\norm{2q_{n_h}\rho}\geq\frac{C}{n_h^2}$, (\ref{e17}), and Proposition \ref{ll8} we get:
\beqano
\parallel(R_{2S_{q_{n_h}}\phi}-Id.)^{-1}\parallel_0\leq Cn_h^{2}.
\eeqano
Moreover, if $\bar{F}=A_h^{(q_{n_{h+1}})}-R_{S_{q_{n_{h+1}}}\phi}$, by Proposition \ref{p2}:
\beqano
\parallel\bar{F}\parallel_0\leq \frac{C}{q_{n_h}}.
\eeqano
In particular, we can apply Proposition \ref{p3} with $\bar{A}:=A^{(q_{n_h})}, \bar{\phi}=S_{q_{n_h}}\phi$ to get 
$B_{h},\tilde{F}\in C^{\infty}(\Bbb{T},SL(2,\Bbb{R})), \tilde{\phi}\in C^{\infty}(\Bbb{T}, \Bbb{R})$ so that \eqref{lp1}--\eqref{lp4} hold. 

By (\ref{e16}) and (\ref{ee18}):
\beqano
\parallel F_h\parallel_1 q_{n_h}\parallel\phi_h\parallel_{50r_0}<1.
\eeqano
Then \eqref{coro1} follows from the estimates of $\parallel F_h\parallel_{50r_0}$ in (\ref{e16}).

\subsection{Going backward from the renormalized cocycle to the starting one} \label{sec.back} 
\nl
Now let $B_h:=\bar{B}_{r_0+1}...\bar{B}_{1}, \tilde{\phi}$ as in Proposition \ref{p4} (with $\bar{A}=A_h^{(q_{n_h})},\bar{\phi}=S_{q_{n_h}}\phi_h$ and with $\bar{B}_1,...,\bar{B}_{r_0+1}$ that are defined applying $r_0+1$ times the cheap trick). We want to show that, because $(\a, A_h)$ commutes with $(q_{n_h}\a, A_h^{(q_{n_h})})$, if $B_h(x+q_{n_h}\a)A_h^{(q_{n_h})}(x)B_h^{-1}(x)$ is close to a rotation valued cocycle, then also $B_h(x+\a)A_h(x)B_h^{-1}(x)$ is close to a rotation valued cocycle. Let also $\bar{B}:=\bar{B}_{r_0+1},\tilde{B}:=\bar{B}_{r_0}...\bar{B}_1$, so that $B_{h}=\bar{B}\tilde{B}$.
\lem{lzl}
\beqano
\bar{B}(x)\tilde{B}(x+q_n\a)\bar{A}(x)\tilde{B}^{-1}(x)\bar{B}^{-1}(x)=R_{\tilde{\phi}(x)}.
\eeqano
\elem{}
\proof
It follows by definition of $\bar{B},\tilde{B}, \phi_h$. \qed
\begin{definition}
Let:
\beqano
J:=
\begin{pmatrix}
0 & -1 \\
1 & 0
\end{pmatrix}
\eeqano
For $M\in SL(2,\Bbb{R})$, 
\beqano
Q(M):=\frac{M+JMJ}{2}.
\eeqano
\end{definition}
In the following Lemma we state some properties of $Q(M)$ that are stated also in \cite{3}.
\lem{az}
For $\theta\in\Bbb{R}, M\in SL(2,\Bbb{R})$, $R_{\theta}Q(M)=Q(R_{\theta}M)=Q(MR_{-\theta})$. Moreover, $M-Q(M)$ is of the form:
\beqano
M-Q(M)=
\begin{pmatrix}
a & b \\
-b & a
\end{pmatrix}.
\eeqano
In particular, $M-Q(M)$ commutes with rotations.
\elem

\nl
\begin{definition}
$\tilde{A}(x):=B_h(x+\a)A_h(x)B_h^{-1}(x), L(x):=Q(\tilde{A}), L_1(x):=Q(\tilde{A}(x+q_{n_h}\a)-\tilde{A}(x))$.
\end{definition}
By Lemma \ref{az}:
\beqano
Q(R_{\tilde{\phi}(x+\a)}\tilde{A}(x)-\tilde{A}(x)R_{\tilde{\phi}(x)})=(R_{\tilde{\phi}(x+\a)}-R_{-\tilde{\phi}(x)})L(x).
\eeqano
\lem{r3}
\beqano
R_{\bar{\phi}(x+\a)}\tilde{A}(x)&-&\tilde{A}(x)R_{\bar{\phi}(x)}=J_1(x)+J_2(x)+J_3(x),
\eeqano
with:
\beqano
J_1(x)&:=&(\bar{B}(x+\a)-\bar{B}(x+\a+q_{n_h}\a))\tilde{B}(x+q_{n_h}\a+\a)A_{h}(x+q_{n_h}\a)B_{h}^{-1}(x+q_{n_h}\a)(R_{\tilde{\phi}(x)}+\tilde{F}(x)),\\
J_2(x)&:=&(\tilde{A}(x+q_{n_h}\a)-\tilde{A}(x))(R_{\tilde{\phi}(x)}+\tilde{F}(x)),\\
J_3(x)&:=&\tilde{A}(x)(\bar{B}(x)-\bar{B}(x+q_{n_h}\a))\tilde{B}(x+q_{n_h}\a)A_h^{(q_{n_h})}(x)B_h^{-1}(x).
\eeqano
\elem
\proof
By Lemma \ref{lzl}
\beqano
R_{\bar{\phi}(x+\a)}\tilde{A}(x)&-&\tilde{A}(x)R_{\bar{\phi}(x)}\\
&=&\bar{B}(x+\a)\tilde{B}(x+q_{n_h}\a+\a)A_h^{(q_{n_h})}(x+\a)A_h(x)B_h^{-1}(x)\\
&-&\bar{B}(x+\a)\tilde{B}(x+\a)A_h(x)\tilde{B}^{-1}(x)\tilde{B}(x+q_{n_h}\a)A_h^{(q_{n_h})}(x)B_h^{-1}(x)\\
&=&\bar{B}(x+\a)\tilde{B}(x+q_{n_h}\a+\a)A_{h}(x+q_{n_h}\a)B_{h}^{-1}(x+q_{n_h}\a)B_{h}(x+q_{n_h}\a)A_{h}^{(q_{n_h})}(x)B_h^{-1}(x)\\
&-&\bar{B}(x+\a)\tilde{B}(x+\a)A_h(x)\tilde{B}^{-1}(x)\tilde{B}(x+q_{n_h}\a)A_h^{(q_{n_h})}(x)B_h^{-1}(x)\\
&=&J_1(x)+J_2(x)+J_3(x).
\eeqano
\qed

\lem{l11l}
Let $t\in\Bbb{N}$. Then:
\beqano
 \parallel(R_{\bar{\phi}(x+\a)}-R_{\bar{\phi}(x)})^{-1}R_{-\bar{\phi}(x)}\parallel_{t}<C(t)n_h^{2(t+1)}\parallel\bar{\phi}\parallel_t.
\eeqano
\elem{}
\proof
It follows by Proposition \ref{ll8} and (\ref{e17}). \qed
\lem{r4}
For $t\leq r_0-7$:
\beqano
\parallel\bar{B}\parallel_{t},\parallel\tilde{B}\parallel_t, \parallel\tilde{\phi}\parallel_t, \parallel\tilde{F}\parallel_t, \parallel\tilde{A}\parallel_{t}, \parallel A_h^{(q_{n_h}}\parallel_t\leq C,
\eeqano
\beqano
\parallel\bar{B}(x)-\bar{B}(x+q_{n_h}\a)\parallel_t\leq C\parallel\tilde{F}\parallel_t.
\eeqano
In particular:
\beqano
\parallel Q(J_1)\parallel_t, \parallel Q(J_3)\parallel_t\leq C\parallel\tilde{F}\parallel_t, \parallel Q(J_2)\parallel_t\leq C\parallel L_1\parallel_t.
\eeqano
\elem
\proof
It follows directly by Proposition \ref{p4}. \qed
\giu
\nl
Let $B_h, \tilde{F}, \tilde{\phi}$ defined as in Proposition 4, $L, L_1, \bar{A}$ as in Definition 4. Let:
\beqano
A_{h+1}(x):=B_{h}(x+\a)A_{h}(x)B_{h}^{-1}(x),\quad R_{\phi_{h+1}}:=\frac{A_{h+1}-L}{(\det(A_{h+1}-L))^{\su{2}}}, F_{h+1}=A_{h+1}-R_{\phi_{h+1}}.
\eeqano
The fact that $R_{\phi_{h+1}}$ is a rotation follows by Lemma \ref{az}.
 The estimates for $B_{h}$ in (\ref{e3}) follow by Proposition \ref{p4}. By Lemma \ref{r3}, Lemma \ref{r4} and Lemma \ref{l11l}, we get the following Lemma:
\lem{10x}
Let $t\leq r_0-7$. Then:
\beqano
\parallel L\parallel_t&\leq& Cn_h^{2(t+1)}\parallel L_1\parallel_t+Cn_h^{2(t+1)}\parallel\tilde{F}\parallel_t\\
&\leq& \frac{Cn_h^{2(t+1)}}{q_{n_h+1}}\parallel L\parallel_{t+1}+Cn_h^{2(t+1)}\parallel\tilde{F}\parallel_t.
\eeqano
\elem 
We also have 
\lem{qw}
For $0\leq{j}\leq r_0-7$:
\beqano
\parallel\tilde{F}\parallel_j\leq\su{q_{n_h+1}^{r_0-j}}.
\eeqano
\elem{}
\proof
It follows by:
\beqano
\parallel\tilde{F}\parallel_j\leq C\parallel\tilde{F}\parallel_{50r_0}^{\frac{j}{50r_0}}\parallel\tilde{F}\parallel_0^{1-\frac{j}{50r_0}},
\eeqano
the estimates of $\parallel\tilde{F}\parallel_0,\parallel\tilde{F}\parallel_{50r_0}$ in (\ref{ee18}), Lemma \ref{lb} and \eqref{coro1}. \qed
From Lemma \ref{10x}, it follows that for $t+l\leq r_0$:
\beq{e19}
\parallel L\parallel_t\leq \frac{Cn_h^{2r_0l(r_0+1)}}{q_{n_h+1}^{l}}\parallel L\parallel_{t+l}+\sum_{j=0}^{l-1}\frac{Cn_h^{2r_0j(r_0+1)}}{q_{n_h+1}^{j}}\parallel\tilde{F}\parallel_{t+j}.
\eeq{}
Moreover:
\beq{e20}
\parallel L\parallel_{r_0}\leq \parallel\tilde{A}\parallel_{r_0}\leq s_{h+1}^6.
\eeq{}
Then, by definition of $F_{h+1}$, Lemma \ref{qw}, (\ref{e19}), (\ref{e20}) it follows that:
\beqano
\parallel F_{h+1}\parallel_{0}<\su{q_{n_{h}+1}^{\frac{r_0}{2}}}, \parallel F_{h+1}\parallel_{50r_0}\leq s_{h+1}^{6}
\eeqano
Finally, by Proposition \ref{p4}, the estimates for $L$ and the definition of $\phi_{h+1}$ and Lemma \ref{lb}, the estimates in Proposition \ref{p1} for $F_{h+1},\phi_{h+1}$ follow.
\qed
\section{Appendix} \label{appendix}
Here we state and prove some basic Lemma that we have used in other sections.
{\bf{Proof of Lemma \ref{lemma1}.}}
Let $t\in\Bbb{N}$, $Y:=\bar{F}R_{-\bar{\phi}}$. Then:
\beqano
|D^tG|_0\leq \sum_{h\geq 1}\left|\frac{D^t Y^h}{h}\right|_0.
\eeqano
By Leibnitz formula, for $h\geq 1$ $D^t Y^h$ is equal to the sum of $h^t$ terms of the form:
\beqano
D^{t_1}Y\dots D^{t_h}Y,
\eeqano
with $t_1+\dots t_h=t$. For each $j$ such that $t_j>0$, by Hadamard's inequality:
\beqano
|D^{t_j}Y|_0\leq C(t)\parallel Y\parallel_0^{1-\frac{t_j}{t}}\parallel Y\parallel_t^{\frac{t_j}{t}}.
\eeqano
So, because there are at most $t$ terms such that $t_j>0$, we get:
\beqano
|D^{t_1}Y\dots D^{t_h}Y|_0\leq C(t)^t\parallel Y\parallel_{0}^{h-1}\parallel Y\parallel_t.
\eeqano
Then:
\beqano
|D^tG|_0\leq C(t) \sum_{h\geq 1}\frac{h^t\parallel Y\parallel_0^{h-1}}{h}\parallel Y\parallel_t \leq C(t)\parallel Y\parallel_t
\eeqano
where in the last inequality we have used that:
\beqano
\parallel Y\parallel_0\leq \parallel\bar{F}\parallel_0\parallel R_{\bar{\phi}}\parallel_0\leq C\parallel\bar{F}\parallel_0\parallel\bar{\phi}\parallel_0\leq \frac{C}{q_n}.
\eeqano
\qed

\lem{v1} For $\bar{A}=R_{\bar{\phi}}+\bar{F}\in C^{\infty}(\Bbb{T}, SL(2,\Bbb{R}))$, let $G$ be such that $\bar{A}=e^G R_{\bar{\phi}}$. Then, for $t\in\Bbb{N}$:
\beqano
|D^t G|_0\leq C(t)\parallel\bar{F}R_{-\bar{\phi}}\parallel_t
\eeqano
\elem{}

\lem{l6}
Let $t\in\Bbb{N}$. There exist polynomials $P_{1,t}(X_1,...X_t),...,P_{t,t}(X_1,...,X_t)$ that are homogenous of degree less or equal then $t$ if the variable $X_i$ has weight $i$ for $i=1,\dots, t$, such that for $g\in C^{\infty}(\Bbb{T},\Bbb{R})$:
\beqano
D^t\left(\su{g}\right)=\sum_{i=1}^{t}\frac{P_{i,t}(Dg,\dots ,D^t g)}{g^{i+1}}.
\eeqano
\elem{}
\lem{v5} {\bf{(See \cite{20})}}
There exists $C>0$ such that, for $0\leq a<b<c\in\Bbb{N}, f\in C^{\infty}(\Bbb{T},\Bbb{R})$:
\beqano
|D^b f|_0\leq C\parallel F\parallel_{a}^{1-\frac{b-a}{c-a}}\parallel F\parallel_{c}^{\frac{b-a}{c-a}}.
\eeqano
\elem{}
\lem{v6}
For $t\in\Bbb{N}$:
\beqano
|D^t R_{\phi}|_0\leq C(t)\parallel\phi\parallel_t.
\eeqano
\elem{}
\proof
It follows by Faa Di Bruno's formula and Hadamard's inequality (Lemma \ref{lemma1}). \qed
\lem{l9}
Let $f\in C^{\infty}(\Bbb{T},\Bbb{R}), a>0, t,h\in\Bbb{N}$. Then:
\beqano
|D^{t+h} T_a(f)|_{0}\leq C(t+h)a^{2+h}\parallel F\parallel_{t},  \quad |D^{t}R_a(f)|_{0}\leq Ca^{-h}\parallel F\parallel_{t+h+2}
\eeqano
\elem{}
\proof
Let $t,h\in\Bbb{N}, a>0$. Then:
\beqano
|D^{t+h}T_a(f)|_{0}\leq \sum_{|l|\leq a} |(2\pi l)|^{t+h}|\hat{f}(l)|\leq |D^{t}f|_0|2\pi a|^{h+2},
\eeqano
where in the last inequality we have used the fact that:
\beqano
|\hat{f}(l)|\leq \frac{|D^{t}f|_0}{|2\pi \bar{l}|^{t}},
\eeqano
with $\bar{l}:=\max\{|l|,1\}$. Now we prove the second inequality:
\beqano
|D^{t}R_a(f)|_{0}\leq \sum_{|l|>a}|\hat{f}(l)\parallel2\pi l|^t\leq \su{a^h}\sum_{|l|>a}|\hat{f}(l)\parallel2\pi l|^{t+h}
\eeqano
\beqano
\leq \frac{|D^{t+h+2}f|_0}{a^h}\sum_{|l|>a}\su{|2\pi l|^2}\leq \frac{C}{a^h}|D^{t+h+2}f|_0.
\eeqano
\qed

\end{document}